\documentclass[letterpaper,11pt]{article}

\usepackage{geometry}
\usepackage[activeacute,english]{babel}
\usepackage{amsmath, amsfonts, amssymb}
\usepackage{amsthm}
\usepackage {graphicx}
\usepackage[all]{xy}

\usepackage{graphics}

\theoremstyle{plain}
 \newtheorem{theorem}{Theorem}[section]
 \newtheorem{lemma}{Lemma}
 \newtheorem{proposition}{Proposition}

\theoremstyle{definition}
 \newtheorem{definition}{Definition}[section]
 
 \newtheorem{example}{Example}[section]
\theoremstyle{remark}
\newtheorem{remark}{Remark}[section]

\newcommand{\bc}{\mathbf{c}}

\newcommand{\mR}{\mathcal{R}}

\newcommand{\bkappa}{\mathbf{\kappa}}

\newcommand{\Z}{\mathbb{Z}}
\newcommand{\R}{\mathbb{R}}
\newcommand{\Q}{\mathbb{Q}}
\newcommand{\rank}{{\rm rank}}

\begin{document}

\title{How far is complex balancing from detailed balancing?}
\author{Alicia Dickenstein and Mercedes P\'erez Mill\'an \\
{\small{Dto. de Matem\'atica, FCEN, Universidad de Buenos Aires,}} \\
{\small{Ciudad Universitaria, Pab. I, C1428EGA Buenos Aires, Argentina.}}\\
 {\small{\texttt{alidick@dm.uba.ar, mpmillan@dm.uba.ar}}}
{ \footnote{Partially supported by UBACYT X064, CONICET PIP 112-200801-00483 and ANPCyT PICT 20569, Argentina}}}

\maketitle

\begin{abstract}
The aim of this article is to build on the use of tools from computational algebra initiated in \cite{TDS}, for the study of general kinetic systems, which have a wide range of applications in chemistry and biology \cite{ctf06,gn08,saf09,so01,gu07,gu08}.
We clarify the relation between the algebraic conditions that must be satisfied by the reaction constants in general (mass action) kinetics systems for the existence of detailed or complex balancing equilibria. The main properties of these systems have been set by Horn, Jackson and Feinberg \cite{fe72,fe89,fe95,ho7273,ho73,ho74,hoja72}. We expect to extend our point of view to the study of qualitative features of the dynamical behaviour of chemical interactions in molecular systems biology.

\smallskip

\end{abstract}

\section{Introduction}

In this article we clarify the relation between the algebraic conditions that must be satisfied by the reaction constants in general mass action kinetics systems, whose main properties have been set by Horn, Jackson and Feinberg \cite{fe72,fe89,fe95,ho7273,ho73,ho74,hoja72}, for the existence of detailed or complex balanced equilibria. These
systems have remarkable dynamic properties and have a wide range of applications in chemistry and biology 
\cite{ctf06,gn08,saf09,so01,gu07,gu08}.

We show that a reversible Horn--Jackson generalized mass action kinetics system satisfying \emph{Feinberg's circuit conditions} is detailed balanced if and only if it is complex balanced. In other words, under \emph{formal} balancing conditions for the cycles (of the underlying undirected graph) of the reaction graph, both notions coincide. We formulate this property in terms of the algebraic equations defining the corresponding varieties in rate constant space.

In order to illustrate some of the definitions and concepts along the paper, we will use a reaction network diagram which represents a nonsequential multisite phosphorylation with two sites. This network models for example the Mek-MKP3-Erk system \cite{gu08,bs97,fb97}:
\begin{small}
    \begin{equation}
\begin{array}{c}
\xymatrix{
 &   &  E + S_{01} \ar@<1ex>[r]^{\kappa_{46}} \ar@<1ex>[dl]^{\kappa_{42}}  & ES_{01} \ar@<1ex>[l]^{\kappa_{64}} \ar@<1ex>[dr]^{\kappa_{67}}  &\\
E + S_{00} \ar@<1ex>[r]^{\kappa_{12}}  & ES_{00} \ar@<1ex>[l]^{\kappa_{21}} \ar@<1ex>[ur]^{\kappa_{24}} \ar@<1ex>[dr]^{\kappa_{23}}  & &   &E + S11 \ar@<1ex>[ul]^{\kappa_{76}} \ar@<1ex>[dl]^{\kappa_{75}}\\
&    & E+S_{10} \ar@<1ex>[r]^{\kappa_{35}} \ar@<1ex>[ul]^{\kappa_{32}}  & ES_{10} \ar@<1ex>[l]^{\kappa_{53}} \ar@<1ex>[ur]^{\kappa_{57}}  &
}\\
\xymatrix{
&    FS_{01} \ar@<1ex>[r]^{\kappa_{1311}}\ar@<1ex>[dl]^{\kappa_{1314}}  & F + S_{01} \ar@<1ex>[l]^{\kappa_{1113}} \ar@<1ex>[dr]^{\kappa_{119}} & & \\
F + S_{00} \ar@<1ex>[ur]^{\kappa_{1413}} \ar@<1ex>[dr]^{\kappa_{1412}} &  & & FS_{11} \ar@<1ex>[r]^{\kappa_{98}}\ar@<1ex>[ul]^{\kappa_{911}}\ar@<1ex>[dl]^{\kappa_{910}} & F + S_{11}\ar@<1ex>[l]^{\kappa_{89}}\\
&  FS_{10} \ar@<1ex>[r]^{\kappa_{1210}}\ar@<1ex>[ul]^{\kappa_{1214}} & F + S_{10} \ar@<1ex>[l]^{\kappa_{1012}} \ar@<1ex>[ur]^{\kappa_{109}} & & 
}
\end{array}
\label{eq:network_1}
\end{equation}
\end{small}
The four phosphoforms, $S_{00}$, $S_{10}$, $S_{01}$, $S_{11}$, are interconverted by the kinase, $E$, and the phosphatase, $F$. Assuming mass action kinetics, each reaction is annotated with the corresponding rate constants, indicated by a choice of numbering of the $14$ complexes in the network. Although $\kappa_{32}$, $\kappa_{42}$, $\kappa_{75}$, $\kappa_{76}$, $\kappa_{109}$, $\kappa_{119}$, $\kappa_{1412}$, $\kappa_{1413}$, are generally taken to be very small and so the corresponding reactions are usually omitted, we will not ignore them in this example because we are interested in special properties of the reaction constants in
reversible networks.

\medskip

In general, a \emph{chemical reaction network} is a finite directed graph whose vertices are labeled by monomials and whose edges are labeled by parameters. The digraph is denoted $G = (V, E)$, with vertex set $V = {1, 2, \dots , n}$ and edge set $E \subseteq \{(i, j) \in V \times V : i \neq j\}$ of cardinal $e$. The node $i$ of $G$ represents the $i$-th chemical complex and is labeled with the monomial
$\bc^{y_i}=c_1^{y_{i1}}c_2^{y_{i2}}\dots c_s^{y_{is}}$. The unknowns $c_1 ,c_2 , \dots , c_s$ represent the concentrations of the $s$ species in the network.
The  $n \times s$-matrix of non-negative integers $Y = (y_{ij} )$ contains the stoichiometric coefficients.
For instance, in the reaction diagram~\eqref{eq:network_1}, the node corresponding to the first complex
$E + S_{00}$ is labeled by the product of the corresponding concentrations $c_E c_{S_{00}}$. In fact, we will name the $12$ concentrations as $c_1, \dots, c_{12}$ 
according to the following order of the species: $S_{00}$, $S_{10}$, $S_{01}$, $S_{11}$, $E$, $F$, $ES_{00}$, $FS_{11}$, $ES_{10}$, $FS_{10}$, $E_{01}$, $FS_{01}$.
So, $c_E c_{S_{00}} = c_5 c_1$. The rows of the stoichiometric matrix $Y \in \{0,1\}^{14 \times 12}$  are ordered according to the numbering of the complexes which is
reflected in the names of the rate constants.

For general systems without a particular (bio)chemical interpretation, 
the name of the species will just be given by a numbering  $\{1, \dots, s\}$ and we record
the monomial labels as  the entries in the row vector 
$$\Psi(\bc) =  (\bc^{y_1} , \bc^{y_2} , \dots , \bc^{y_n}).$$
In our example, $\Psi(\bc) =  (c_5c_1, c_7, c_5c_2, c_5c_3, c_9, c_{11}, c_5c_4, c_6c_4, c_8, c_6c_2, c_6c_3, c_{10}, c_{12}, c_6c_1).$

Each directed edge $(i,j) \in E$ is labeled by a positive parameter $\kappa_{ij}$ which represents the rate constant in the reaction from the $i$-th chemical complex to the $j$-th chemical complex. Note that if there is an edge from $i$ to $j$ and an edge from $j$ to $i$ then we have two unknowns $\kappa_{ij}$ and $\kappa_{ji}$.

Such system will be denoted by $G=(V,E,\bkappa,Y)$, with $\bkappa=(\kappa_{ij})_{(i,j) \in E}$ the vector in $\R_{>0}^{E}$ of the rate constants. By abuse of notation, we will just write a chemical reaction network as a digraph $G=(V,E)$. Also, we will assume that a numbering of $E$ is chosen and we will consider $\bkappa \in \R_{>0}^{e}$.

Let $A_\kappa$ denote the negative of the \emph{Laplacian} of $G$. Hence $A_\kappa$ is the $n \times n$-matrix whose off-diagonal entries are the $\kappa_{ij}$ and whose row sums are zero. Mass action kinetics specified by the digraph $G$ is the dynamical system
\begin{equation} \label{dynsys}
 \frac{d{\bc}}{dt}=\Psi(\bc) A_\kappa Y,
\end{equation}
where $\bc(t) = (c_1(t), \dots, c_s(t))$.

\medskip

\begin{definition}
 A \emph{complex balanced} system is a dynamical system ~\eqref{dynsys} for which the algebraic equations $\Psi(\bc)A_\kappa = 0$ admit a strictly positive solution $\bc_0 \in \R_{>0}^s$ . Such a solution $\bc_0$ is a steady state of the system, i.e., the $s$ coordinates of $\Psi(\bc_0)A_\kappa Y$ vanish.
\end{definition}

\begin{remark}
Clearly, a system ~\eqref{dynsys} being complex balanced depend on both the digraph $G$ and the rate constants $\kappa_{ij}$. A main property of complex balanced systems is that \emph{all} strictly positive steady states $c$
satisfy $\Psi(\bc)A_\kappa = 0$.
They are \emph{quasi-thermostatic} \cite{hoja72}, which in the terminology of \cite{TDS} means that the positive steady state variety is \emph{toric}.
\end{remark}


From now on we will assume that the digraph $G$ is reversible, i.e. if $(i,j) \in E$, then $(j,i) \in E$. We can thus identify $G = (V, E)$ with the underlying undirected graph $\widetilde{G} = (V, \widetilde{E})$, where $\widetilde{E} = \{\{i, j\} : (i, j) \in E\}$.

\begin{definition}
A detailed balanced system is a dynamical system ~\eqref{dynsys} for which the following algebraic equations admit a strictly positive solution $\bc_0 \in \R_{>0}^s$:
\begin{equation} \label{eq:dbal}
- \kappa_{ij} \bc_0^{y_i} + \kappa_{ji} \bc_0^{y_j} = 0, \quad {\text{for \ all \, }} \, \{i, j\} \in \widetilde{E}.
\end{equation}
As it is for complex balanced systems, the condition of being detailed balanced depends on the graph $\widetilde{G}$ and the constants $\kappa_{ij}$.
\end{definition}

Note that $A_\kappa$ decomposes as the sum of $n \times n$
matrices $A_\kappa^{\{ij\}}$  for each undirected edge $\{i, j\} \in \widetilde{E}$ of the graph $G$, where in rows $i, j$ and columns $i, j$ the matrix $A_\kappa^{\{i,j\}}$ equals
$$\left(\begin{array}{cc}
-\kappa_{ij} & \kappa_{ij}\\
 \kappa_{ji} & -\kappa_{ji}
 \end{array} \right),$$
and all other entries of the matrix $A_\kappa^{\{i,j\}}$ are $0$. Since the algebraic equation $- \kappa_{ij} \bc_0^{y_i} + \kappa_{ji} \bc_0^{y_j} = 0$ means
that $\Psi(\bc_0) A_\kappa^{\{i,j\}} = 0$,  we see that every detailed balanced system is also complex balanced.
The converse is not true in general. Again, a main
property of a detailed balanced system is that \emph{all} of its positive steady states $\bc$ satisfy
$- \kappa_{ij} \bc^{y_i} + \kappa_{ji} \bc^{y_j} = 0$ for $\{i, j\} \in \widetilde{E}$. 

\medskip

For instance, in our example~\eqref{eq:network_1}, for {\em any} choice of $\kappa \in \R_{>0}$ the system with the following rate constants 
is complex balanced but not detailed balanced:
\begin{small}
\begin{equation}\label{eq:constants}
 \begin{array}{l}
  \kappa_{12}=\kappa_{24}=\kappa_{46}=\kappa_{53}=\kappa_{67}=\kappa_{89}=\kappa_{910}=\kappa_{1012}=\kappa_{1214}=\kappa_{1311}=\kappa\\
\\
  \kappa_{32}=\kappa_{42}=\kappa_{75}=\kappa_{76}=\kappa_{109}=\kappa_{119}=\kappa_{1412}=\kappa_{1413}=\frac{1}{4}\kappa\\
\\
  \kappa_{23}=\kappa_{35}=\kappa_{57}=\kappa_{64}=\kappa_{911}=\kappa_{1113}=\kappa_{1314}=\frac{3}{4}\kappa, \, \,
  \kappa_{21}=\frac{23}{4}\kappa, \, \,
  \kappa_{98}=\frac{47}{4}\kappa, \, \,
  \kappa_{1210}=\frac{69}{22}\kappa.
 \end{array}
\end{equation}
\end{small}
\noindent For any $\alpha \in \mathbb{R}_{>0}$,  $\bc_{0, \alpha}=\alpha \left(23,17,11,47,1,2,4,8,14,11,13,16\right)$ is a positive steady state of the system for  which $\Psi(\bc_{0, \alpha})A_\kappa=0$, and hence the system is complex balanced. On the other side, for this choice of rate constants the system is {\em not} detailed balanced since~\eqref{eq:dbal} does not hold for instance for $i=4, \, j=6$, that is, for the pair of reactions
$$\begin{array}{c}
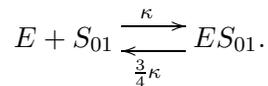
\xymatrix{
 E + S_{01} \ar@<1ex>[r]^{\kappa}  & ES_{01}. \ar@<1ex>[l]^{\frac{3}{4}\kappa}  \\}
\end{array}$$

\medskip

We now recall Feinberg's circuit conditions \cite{fe89}, also known as Wegscheider's condition \cite{ss89}.
Since they correspond to linear relations which only depend on the structure of the reaction graph and not
on the particular complexes, we prefer to name them as \emph{formal} conditions.
For every cycle $\widetilde{C}$ in $\widetilde{G}$, we will choose one direction and define $C^+$ as the cycle in $G$ in that direction. $C^-$ will be the cycle in the opposite direction. Although the directions are arbitrarily chosen, we will not worry about that since we will only need to distinguish between the two of them.


\begin{definition} \label{def:formal}
A formally balanced system is a dynamical system ~\eqref{dynsys} for which the following condition holds for every cycle $\widetilde{C}$ of $\widetilde{G}$:
\begin{equation}\label{eq:fbq}
 \underset{(i,j) \; in \; C^+}{\prod}\kappa_{ij}=\underset{(j,i) \; in \; C^-}{\prod}\kappa_{ji}.
\end{equation}
\end{definition}

We will talk about formally balanced systems although this definition can be applied to any digraph whose edges are reversible and labeled by constants $\kappa_{ij}$.

In our example~\eqref{eq:network_1}, we can consider the cycle $\widetilde{C}$:
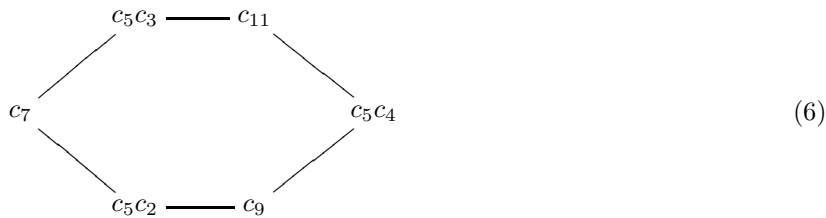
\begin{small}
    \begin{equation}
\begin{array}{c}
\xymatrix{
   &  c_5c_3 \ar@{-}[r] \ar@{-}[dl]  & c_{11} \ar@{-}[l] \ar@{-}[dr]  &\\
  c_7  & &   &c_5c_4 \\
   & c_5c_2 \ar@{-}[r] \ar@{-}[ul]  & c_9 \ar@{-}[l] \ar@{-}[ur]  &
}
\end{array}
\label{eq:cycle}
\end{equation}
\end{small}
As~\eqref{eq:fbq} is not satisfied for $\widetilde{C}$, the system is not formally balanced.

\medskip

Equations ~\eqref{eq:fbq} show that the set
\[\mathcal{FB}_Y=\{\bkappa=(\kappa_{ij})_{(i,j) \in E}: G=(V,E,\bkappa,Y) \; is \; formally \; balanced\}\]
is an algebraic variety in $\R_{>0}^{e}$.

It follows from \cite{fe89}  that the set
\[\mathcal{DB}_Y=\{\bkappa=(\kappa_{ij})_{(i,j) \in E}: G=(V,E,\bkappa,Y) \; is \; detailed \; balanced\}\]
is also an algebraic variety in $\R_{>0}^{e}$ (see Lemma~\ref{lem:db}).

In turn, it follows from \cite[Section 2]{TDS} that the set
\[\mathcal{CB}_Y=\{\bkappa=(\kappa_{ij})_{(i,j) \in E}: G=(V,E,\bkappa,Y) \; is \; complex \; balanced\}\]
is a third algebraic subvariety of  $\R_{>0}^{e}$ (see Proposition~\ref{prop:cb}), called the \emph{moduli space of toric dynamical systems} in \cite{TDS}.

As we have already remarked, $ \mathcal{DB}_Y \subseteq \mathcal{CB}_Y$. In fact, the main Theorem in
\cite{fe89} shows that  $ \mathcal{DB}_Y \subseteq \mathcal{FB}_Y$.

\medskip

In this paper we prove  the following result for a mass action kinetics dynamical system associated to
a reversible
chemical reaction network $G=(V,E, \bkappa,Y)$:

\begin{theorem} \label{th:main}
Under the assumption of formal balancing, a reversible system is detailed balanced
if and only if it is complex balanced. That is,
\begin{equation} \label{eq:main}
\mathcal{CB}_Y \cap \mathcal{FB}_Y = \mathcal{DB}_Y.
\end{equation}
\end{theorem}

Our result generalizes two particular situations in which it is known
that the notions of detailed and complex balancing coincide:
the case in which $\widetilde{G}$ has no
cycles, and the case of deficiency zero networks for which
 $\mathcal{DB}_Y = \mathcal{FB}_Y$ (\cite{fe89}, see also Lemma~\ref{lem:fb} below).

Our algebraic approach follows the lines of \cite{TDS}.
Our arguments easily imply that \eqref{eq:main} holds at the level of ideals (which are radical).
We refrain from giving a more algebraic formulation since it is straightforward and our main concern is
to clarify these notions in the framework of general mass action kinetics networks.
In Section~\ref{sec:prelim} we recall known results, mainly from \cite{fe89,TDS},
that we state in a language adapted to our setting.  In particular,
we introduce new quotient variables which  allow us to organize the proof of
Theorem~\ref{th:main} in  Section~\ref{sec:main}. 

In Section~\ref{sec:general}, following a suggestion of Martin Feinberg,  we translate Theorem~\ref{th:main} to the setting of general kinetic systems in Theorem~\ref{th:maingeneral} and we express another necessary and
sufficient condition for a complex balanced system to be detailed balanced in Proposition~\ref{prop:one}.

\section{Preliminaries} \label{sec:prelim}

Given a chemical reaction network $G=(V,E,\bkappa,Y)$,  we will denote by
$\widetilde{G}$ the associated undirected graph. As we assume that $G$ is reversible,
there is no loss of information in passing to $\widetilde{G}$.

\subsection{The minors of $A_\kappa$}\label{ssec:kappa}

Let $G$ be a reversible digraph corresponding to a chemical reaction network
and call  $G_t=(V_t,E_t)$, $t=1,\dots,\ell,$ the connected components of $G$.
Up to renumbering, we can assume  $A_\kappa = A_\kappa(G)$ is block diagonal, with diagonal blocks
the corresponding matrices $A_\kappa(G_t)$ for the components $G_1, \dots, G_\ell$.
Following \cite{TDS}, we introduce the following definition:

\begin{definition} \label{def:Ki}
Consider any directed subgraph $T= (V(T), E(T))$ of $G$ with $n - \ell$ edges whose
underlying undirected graph
is a spanning forest of the the underlying undirected graph  of $G$.
Fix a connected component $G_t$ of $G$ and
write $\kappa^T_t$ for the product of the $\# V_t -1$ rate constants which correspond to
all edge labels of the edges in $E(T) \cap E_t$.
Let $i$ be one of the nodes of $G_t$. The directed tree obtained by the restriction $T_t$ of $T$ to $G_t$
is called an $i$-tree if the node $i$ is its unique sink, i.e., all edges
are directed towards node $i$.
We will write $\kappa^{T_t}$ for the product of the $\#V_t-1$ rate constants which
correspond to all edge labels of the edges of $T_t$. We introduce the following
constants, which are polynomials in the $(\kappa_{ij})$:

\begin{equation} \label{eq:Ki}
K_i=\underset{T_t\; an \; i-tree}{\sum}\kappa^{T_t}.
\end{equation}
Note that each $K_i$ is a nonempty  sum of positive terms, because as $G_t$ is strongly connected,
there exists at least one $i$-tree for every vertex $i$ and each $\kappa_{uv}>0$ for $(u,v) \in E_t$.
\end{definition}

It follows from the Matrix-Tree Theorem \cite{st99} that for any $i \in V_t$, the absolute value of
the determinant of the submatrix of $A_\kappa(G_t)$ obtained by deleting the $i$-th row and any one of
the columns, equals $K_i$. This (non-zero) minor is independent (up to sign) of the choice of columns
because the row sums of $A_\kappa(G_t)$ are zero. Compare also with the statements in \cite{ka56}.

\begin{example}\label{infinito}
We will introduce a new mathematical example only for the purpose of making the calculations more transparent. Let $G=(\{1,2,3,4,5,6\},E,\bkappa,Y)$ be the following connected chemical reaction network:
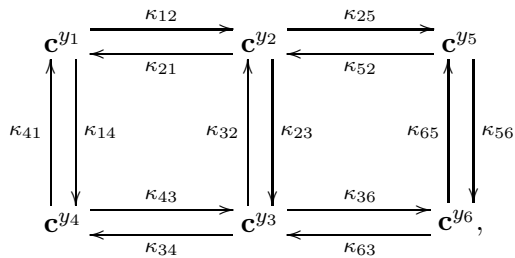
\begin{equation*}
\xymatrix{
\bc^{y_1} \ar@<1ex>[dd]^{\kappa_{14}} \ar@<1ex>[rr]^{\kappa_{12}} & & \bc^{y_2} \ar@<1ex>[ll]^{\kappa_{21}} \ar@<1ex>[dd]^{\kappa_{23}} \ar@<1ex>[rr]^{\kappa_{25}} & & \bc^{y_5} \ar@<1ex>[ll]^{\kappa_{52}} \ar@<1ex>[dd]^{\kappa_{56}}\\
\\
\bc^{y_4} \ar@<1ex>[uu]^{\kappa_{41}} \ar@<1ex>[rr]^{\kappa_{43}} & &  \bc^{y_3} \ar@<1ex>[ll]^{\kappa_{34}} \ar@<1ex>[uu]^{\kappa_{32}} \ar@<1ex>[rr]^{\kappa_{36}} & & \bc^{y_6} \ar@<1ex>[ll]^{\kappa_{63}} \ar@<1ex>[uu]^{\kappa_{65}},
}
\end{equation*}
with
$E=\{(1,2),(1,4),(2,3),(2,5),(3,4),(3,6),(5,6),\-(2,1),(4,1),(3,2),(5,2),(4,3),(6,3),(6,5)\}$.

For example, $K_1=\underset{T\; an \; 1-tree}{\sum}\kappa^T= \kappa_{21}\kappa_{32}\kappa_{63}\kappa_{41}\kappa_{52}+\kappa_{21}\kappa_{32}\kappa_{63}\kappa_{41}\kappa_{56}+\kappa_{21}\kappa_{32}\kappa_{63}\kappa_{43}\kappa_{52}$

$+\kappa_{21}\kappa_{32}\kappa_{63}\kappa_{43}\kappa_{56}+\kappa_{21}\kappa_{52}\kappa_{65}\kappa_{32}\kappa_{41}+\kappa_{21}\kappa_{52}\kappa_{65}\kappa_{32}\kappa_{43}+\kappa_{21}\kappa_{52}\kappa_{65}\kappa_{34}\kappa_{41}+\kappa_{21}\kappa_{52}\kappa_{65}\kappa_{36}\kappa_{41}$

$+\kappa_{21}\kappa_{52}\kappa_{65}\kappa_{36}\kappa_{43}+\kappa_{41}\kappa_{63}\kappa_{34}\kappa_{21}\kappa_{52}+\kappa_{41}\kappa_{63}\kappa_{34}\kappa_{21}\kappa_{56}+\kappa_{41}\kappa_{63}\kappa_{34}\kappa_{23}\kappa_{52}+\kappa_{41}\kappa_{63}\kappa_{34}\kappa_{23}\kappa_{56}$

$+\kappa_{41}\kappa_{63}\kappa_{34}\kappa_{25}\kappa_{56}+\kappa_{41}\kappa_{52}\kappa_{23}\kappa_{34}\kappa_{65} $.
\end{example}

\subsection{The linear relations} \label{ssec:la}

We now recall some results from linear algebra, which in a different language are all contained in
\cite{fe89}.
Choose a numbering of the set of reactions, that is of the set of edges $E$ of $G$ and form
the signed incidence matrix $C_G \in \{-1, 0, 1\}^{n \times e}$ whose column associated to the reaction
$(i,j)$ has a $-1$ on row $i$, a $1$ on row $j$ and all other entries equal to $0$. We denote by
$$ N \, = \, \ker_\Z (Y^t \cdot C_G) \, \subset \Z^e, $$
the kernel over $\Z$ of the product matrix $Y^t \cdot C_G$.

Clearly, the vector with a $1$ on its $(i,j)$-th coordinate and on its $(j,i)$-th coordinate
lies in $N$.  We could instead choose one direction for each pair of reactions $(i,j), (j,i)$ in any way to
get a directed subgraph $G'$ of $G$,
and consider the associated signed incidence matrix $C_{G'}$, with integer kernel $N' \subset \Z^{\frac e 2}$,
since we can clearly reconstruct $N'$ from $N$ and vice versa.

The following combinatorial arguments go back to Kirchoff. We can distinguish the following
sublattice $N'_1$ of $N'$. It is the $\Z$-module spanned by the cycles of the underlying
undirected graph $\widetilde{G}$. More precisely, given any oriented cycle $\mathcal C$
we form the vector $v_{\mathcal C} \in \{-1, 0, 1\}^{\frac e 2}$ whose $(i,j)$ coordinate equals
$1$ if  the edge $(i,j) \in G'$ is in $\mathcal C$, $-1$ if instead the edge $(j, i)$ lies in $\mathcal C$,
and $0$ if neither of the edges $(i,j), (j,i)$ is in $\mathcal C$.

The rank of $N'_1$ equals
$ \frac e 2- n +\ell$, and a basis is formed by the \emph{fundamental cycles} associated to a choice of
a spanning forest $T$ of $G$. The fundamental cycles associated to $T$ are those (undirected) cycles
which are created when we add an edge in the associated undirected graph $\widetilde{T}$
between any two vertices in the same connected component of $G$. Note that
although the number of fundamental cycles in a graph is fixed,
the cycles that become fundamental change with the spanning forest.

If we fix a spanning forest $\widetilde{T}$ of $\widetilde{G}$, we can moreover
choose a direct complement $N'_2$ of $N'_1$ in $N'$ as follows.
Consider all vectors $v=(v_{ij}, \, (i,j) \in E(G'))$ in $N'$ such that
$v_{ij} \not= 0 \Rightarrow \{i,j\} \in E(\widetilde{T})$. Call
$N'_2$ the $\Z$-span of all these vectors $v$ with support contained in $E(\widetilde{T})$.
Then
\[ N' \, = \, N'_1 \oplus N'_2.\]
The concept of deficiency $\delta$ of a chemical reaction network has been introduced and studied by
M. Feinberg in a series of papers \cite{fe79,fe89}. With our notations, the deficiency of
the network $G$ equals $\delta = n - \dim S - \ell$, where $S$ is the stoichiometric linear space defined
by
\[ S \, = \, {\rm span} \{ y_i - y_j, \, (i,j) \in E \}.\]
Thus, $\dim S = \rank (Y^t \cdot C_{G'}) = \rank (Y^t \cdot C_G)$. As $\dim S = \frac e 2 - \rank (N')$, we get that
 $\rank (N'_2) = \delta$, so that $N'_2=0$ if and only if $\delta =0$, and for $\delta >0$
 we could choose a system of $\delta$ generators of $N'_2$.

In a similar way, we can decompose $N$ as $N = N_0 \oplus N_1 \oplus N_2$, where $N_0$ is the lattice of
rank $\frac e 2$ spanned by the $0,1$ vectors in $N$ which express the fact that the $(i,j)$-th column
of $Y^t \cdot C_G$ is minus its $(j,i)$-th column, and $N_i$ for $i=1,2$ is isomorphic to
$N'_i$ (we simply add $0$ coordinates for the entries corresponding to the edges not in $G'$).

\subsection{Complex and detailed balanced systems}\label{ssec:cdb}

We will now present some equivalencies for complex and detailed balanced systems. We will use the following basic
result \cite{es96}:

\begin{lemma} \label{lem:binom}
Let $\bf k$ be a field and $a_1, \dots, a_m \in \Z^s$. Given a vector $z=(z_1, \dots, z_m) \in ({\bf k} - \{0\})^m$,
there exists $x=(x_1, \dots, x_s) \in ({\bf k} - \{0\})^s$ such that $z_i = x^{a_i}$ for all $i=1, \dots,m$ if
and only if $z^\lambda = 1$ for all $\lambda \in \Z^m$ such that $\sum_{i=1}^m \lambda_i a_i =0$.
\end{lemma}
Here, we use the standard notation $z^\lambda = \prod_{i=1}^m z_i^{\lambda_i}$.
When ${\bf k} = \R$, an easy proof can be given by taking logarithms.

We now introduce new variables which are suitable for our formulation.

\begin{definition}
Let $G=(V,E,\bkappa,Y)$ be a reversible chemical reaction network defining a dynamical system as in \eqref{dynsys}.
For each $(i,j) \in E$ we define
\begin{equation}\label{eq:qK}
q_{ij} \, = \, \frac{\kappa_{ij}}{\kappa_{ji}}, \quad Q_{ij} \, = \, \frac{K_j}{K_i}.
\end{equation}
\end{definition}

\begin{remark} \label{rmk:qij}
The following equations hold
\[ q_{ij} q_{ji} \, = \, Q_{ij} Q_{ji} \, = 1, \quad {\rm \, for \, all \, } \, (i,j) \in E.\]
We define $Q_{ij}$ by the same formula for any pair $i,j$ in $1, \dots, n$ and then
\[ Q_{ij} Q_{jk} \, = Q_{ik}, \quad {\rm \, for \ all \, } \, i,j,k \in \{1, \dots, n\}.\]
\end{remark}

We easily have

\begin{lemma} \label{lem:db}
A chemical reaction network, $G=(V,E,\bkappa,Y)$ is detailed balanced if and only if
\begin{equation}\label{eq:dbq}
q^\lambda = 1, \quad {\rm \, for \ all \, }\,  \lambda \in N.
\end{equation}
Here $q$ denotes the vector $q = (q_{ij})_{(i,j) \in E}$.
\end{lemma}

\begin{proof}
Clearly, a positive vector $c_1$ satisfies a binomial equation
$-\kappa_{ij} \bc_0^{y_i} + \kappa_{ji} \bc_0^{y_j}=0$ if and only
if $\bc_0^{y_i-y_j} = q_{ij}$.
The result follows from
Lemma~\ref{lem:binom} for $m=e$ and $\{a_1, \dots, a_m\} = \{ y_i - y_j, (i,j) \in E\}$.
\end{proof}

We can also translate in a straightforward fashion our definition of formal balancing.

\begin{lemma} \label{lem:fb}
Given a chemical reaction network, $G=(V,E,\bkappa,Y)$,
the following statements are equivalent:
\begin{itemize}
\item[(i)] The associated system is formally balanced,
\item[(ii)] For every cycle $\widetilde{C}$ of $\widetilde{G}$ it holds that
\begin{equation}\label{eq:fbq2}
 \underset{(i,j) \; in \; C^+}{\prod} q_{ij} \, = \, 1,
\end{equation}
\item[(iii)] The vector $q = (q_{ij})_{(i,j) \in E}$ verifies
\begin{equation}\label{eq:fbq3}
q^\lambda = 1, \quad {\rm \, for \ all \, }\,  \lambda \in N_1.
\end{equation}
\end{itemize}
\end{lemma}

Then, a formally balanced system $G$ is detailed balanced if and only if
Equations \eqref{eq:dbq} hold for all $\lambda$ in a set of generators of $N_2$.

Using the results in \cite[Section 2]{TDS} we also have:

\begin{proposition} \label{prop:cb}
A chemical reaction network, $G=(V,E,\bkappa,Y)$ is complex balanced if and only if
\begin{equation}\label{eq:cbq}
Q^\lambda = 1, \quad {\rm \, for \ all \, }\,  \lambda \in N_2.
\end{equation}
Here, $Q$ denotes the vector $Q = (Q_{ij})_{(i,j) \in E}$.
\end{proposition}

\begin{proof}
We first claim that a network $G= (V, E,\bkappa, Y)$ defines a complex balanced system
if and only if there exists a positive vector $\bc_0 \in \R^s$ such that the following
binomial equations are satisfied
\begin{equation} \label{eq:K}
K_i\bc_0^{y_j}-K_j\bc_0^{y_i}=0, \, {\rm \, for \, all\, } \,  (i,j) \in E.
\end{equation}
It is easy to see, as in Remark~\ref{rmk:qij}, that Equations \eqref{eq:K} are
equivalent to
\begin{equation}\label{eq:Kbis}
 K_i\bc_0^{y_j}=K_j\bc_0^{y_i}, \, {\rm \, for \, all\, } \,  (i,j) \in G_t, {\rm \, for \,
some \,} \, t = 1, \dots,\ell.
\end{equation}

To prove that complex balancing is equivalent to \eqref{eq:Kbis},
 we form as in \cite{TDS} the following binomial ideals in
$\Q[\bc]:=\Q[c_1, \dots, c_s]$:
\begin{equation}\label{eq:I}
I = I_1 + \dots + I_t, \quad I_t = \langle  K_i c^{y_j} - K_j c^{y_i}, (i,j) \in E_t \rangle, \, t=1, \dots, \ell.
\end{equation}
We moreover define the ideals $T_G$ as the saturation
\[T_G=(I:(c_1c_2\dots c_s)^\infty)=\{p \in \Q[\bc]: \exists \; u \in \Z_{\geq 0} {\rm \; such \; that \;} p(c_1c_2\dots c_s)^u \in I\}.\]
We denote by  $V_{>0}(I)$ the positive variety of $I$, that is, the zeros of $I$ in $(\R_{>0})^s$, and similarly
for other ideals.
As $ T_G=(I:(c_1c_2\dots c_s)^\infty) = (I_1:(c_1c_2\dots c_s)^\infty)+ \dots
+ (I_\ell:(c_1c_2\dots c_s)^\infty),$
we deduce from display (8) in \cite{TDS} that
$V_{>0}(T_G)=\{\bc \in \R_{>0}^{s} : \Psi(\bc)A_\kappa=0\}$.
But a point $x$ with all non zero coordinates is annihilated by $T_G$ if and only if it is annihilated by $I$.
We then have that
\[V_{>0}(I)=\{\bc \in \R_{>0}^{s} : \Psi(\bc)A_\kappa=0\}, \]
and so the system $G$ is complex balanced if and only if there exists a positive vector $\bc_0$ satisfying
equations ~\eqref{eq:K}.

Now, we argue as in the proof of Lemma~\ref{lem:db}. These equations are equivalent to
$c_0^{y_i-y_j} = Q_{ij}$, for all $(i,j) \in E$. By Lemma~\ref{lem:binom} for $m=e$ and $\{a_1, \dots, a_m\} = \{ y_i - y_j, (i,j) \in E\}$,
these conditions are in turn equivalent to
$$Q^\lambda = 1, \quad {\rm \, for \ all \, }\,  \lambda \in N.$$
But from the definition of the vector $Q$, it is clear that the equalities $Q^\lambda = 1$ always hold for any $\lambda \in N_0 \cup N_1$, 
and so the result follows.
\end{proof}

\begin{remark}
Note that in fact it is enough to check equalities~\eqref{eq:cbq} for $\lambda$ in a basis of $N_2$.
For instance, the rank of $N_2$   in  the network ~\eqref{eq:network_1} is $3$. It is straightforward to check that
for any choice of constants as in~\eqref{eq:constants}:
$$Q_{12}^1 \times Q_{24}^1 \times Q_{1113}^1 \times Q_{1314}^1=\frac{K_1K_{11}}{K_4K_{14}}=1$$
$$Q_{12}^1 \times Q_{23}^1 \times Q_{1012}^1 \times Q_{1214}^1=\frac{K_1K_{10}}{K_3K_{14}}=1$$
$$Q_{35}^1 \times Q_{57}^1 \times Q_{89}^1 \times Q_{910}^1=\frac{K_3K_8}{K_7K_10}=1,$$
which proves again that the system is complex balancing (without needing to show a complex balancing
steady state).
\end{remark}

\section{Proof of Theorem~\ref{th:main}} \label{sec:main}

Consider a formally balanced
reversible chemical reaction network $G = (V,E,\bkappa,Y)$.
By Lemmas~\ref{lem:db}  and~\ref{lem:fb} and Proposition~\ref{prop:cb}, we need to show that
if the constants $q_{ij}$ satisfy Equations \eqref{eq:fbq2}, then
\[ Q^\lambda =1, \quad {\rm \, for \ all \, }\,  \lambda \in N.\]
if and only if
\[ q^\lambda =1, \quad {\rm \, for \ all \, }\,  \lambda \in N.\]

These relations possibly involve constants associated to edges in several
connected components of $G$.
In fact, it holds that modulo the formal balancing relations,
an algebraic dependency relation $P(K) =0$ among the (invertible) variables
$Q_{ij}$ holds for a polynomial $P$ in $e$ variables
if and only if the ``same'' algebraic relation $P(q)=0$
is true for the (invertible) variables $q_{ij}$. This follows immediately from
the following Proposition.

\begin{proposition} \label{prop:key}
 Let $G = (V,E,\bkappa,Y)$
be a  reversible mass action kinetics system which is formally balanced. Then,
\begin{equation}\label{eq:key}
 Q_{ij} \, = \, q_{ij}, \quad {\rm \, for \, all \,} \, (i,j) \in E.
\end{equation}
\end{proposition}

\medskip

\begin{proof}
Since  Equations~\eqref{eq:fbq2} relate variables $q_{uv}$ for $(u,v)$ in a single
connected component of $G$, and given $(i,j) \in E$ then $i,j$ belong to the same component,
we can assume $G$ is connected.

Fix $(i,j) \in E$. We define a bijection between the set of
$j$-trees and the set of $i$-trees as follows (see Example~\ref{ex:continued}
for an illustration).
Let $T$ be any $j$-tree.

\begin{itemize}
\item[(i)] If the edge $(i,j) \in E(T)$, then let $T'$
be the tree obtained by replacing $(i,j)$ by the opposite edge $(j,i)$.
\item[(ii)] If the edge $(i,j) \notin E(T)$, let ${\mathcal C}_{ij}$ be the \emph{undirected}
fundamental cycle which is created in $\widetilde{T}$ by adding the edge $(i,j)$.
Call ${\mathcal C}_{ij}^+$ the corresponding \emph{oriented} cycle which contains $(i,j)$.
Then, let $T'$ be the tree obtained by  giving to the edges
of $T$ which ``lie'' on ${\mathcal C}_{ij}$ the direction in ${\mathcal C}_{ij}^+$
(that is, we ``reverse'' all these edges in $T$).
\end{itemize}

It is straightforward to check that in both cases $T'$ is in fact an $i$-tree and that
the map $T \mapsto T'$ is a bijection. So,
we have established a bijection between the terms in $K_i$ and the terms in $K_j$.

Let $T$ be a $j$-tree. We now compare the term $\kappa^T$ in $K_j$
with the corresponding term $\kappa^{T'}$ in $K_i$. If $(i,j) \in E(T)$,
we clearly have that
\[ \kappa^T \, = \, q_{ij} \, \kappa^{T'}.\]
If instead we have that $(i,j) \notin T$ then
\[  \kappa^T = \left(\prod_{(u,v) \in {\mathcal C}_{ij}^+, (u,v) \not= (i,j)} q_{vu}\right) \, \kappa^{T'}.\]
Now, by the assumption of formal balancing, we have that $\underset{(u,v) \in {\mathcal C}_{ij}^+}{\prod} q_{uv} =1$
and so
\[\prod_{(u,v) \in {\mathcal C}_{ij}^+, (u,v) \not= (i,j)} q_{vu} = q_{ij}.\]
Therefore,
\[ Q_{ij} \, = \frac{K_j} {K_i} \, = \, q_{ij},\]
as wanted.
\end{proof}

\begin{example}[Example~\ref{infinito} continued] \label{ex:continued}
Choose the following $1$-tree $T$:
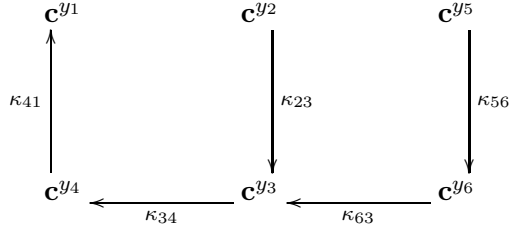
$$\xymatrix{
\bc^{y_1}   & & \bc^{y_2} \ar@<1ex>[dd]^{\kappa_{23}}  & & \bc^{y_5}  \ar@<1ex>[dd]^{\kappa_{56}}\\
\\
\bc^{y_4} \ar@<1ex>[uu]^{\kappa_{41}}  & &  \bc^{y_3} \ar@<1ex>[ll]^{\kappa_{34}}   & & \bc^{y_6} \ar@<1ex>[ll]^{\kappa_{63}}
}$$

Let $(i, j)=(4,1)$. It is clear that by reversing the edge $(4,1) \in E(T)$ one gets a $4$-tree.

Let now $(i, j)=(2,1)$, which does not lie in $E(T)$,   and ${\mathcal C}_{12}^+$
be the corresponding oriented fundamental cycle:
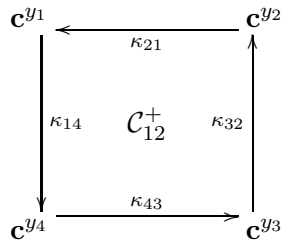
\begin{center}
$\begin{array}{c}
\xymatrix{
\bc^{y_1} \ar@<1ex>[dd]^{\kappa_{14}}  & & \bc^{y_2} \ar@<1ex>[ll]^{\kappa_{21}}  \\
& {\mathcal C}_{12}^+ & \\
\bc^{y_4}  \ar@<1ex>[rr]^{\kappa_{43}} & &  \bc^{y_3}  \ar@<1ex>[uu]^{\kappa_{32}}
}
\end{array}$
\end{center}
Then, reversing the arrows in the cycle gives the following $2$-tree $T'$
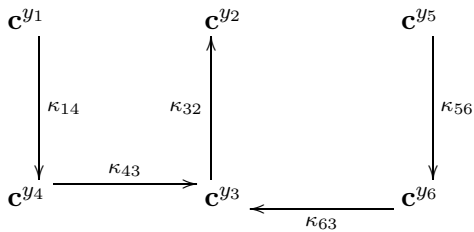
$$\xymatrix{
\bc^{y_1} \ar@<1ex>[dd]^{\kappa_{14}}  & & \bc^{y_2}   & & \bc^{y_5}  \ar@<1ex>[dd]^{\kappa_{56}}\\
\\
\bc^{y_4} \ar@<1ex>[rr]^{\kappa_{43}} & &  \bc^{y_3} \ar@<1ex>[uu]^{\kappa_{32}}   & & \bc^{y_6} \ar@<1ex>[ll]^{\kappa_{63}}
}$$
\end{example}

\section{General kinetic systems} \label{sec:general}

In this section we generalize Theorem \ref{th:main} to non necessarily mass-action kinetic systems in the sense of \cite{fe79}, see also \cite[Section 2]{saf09}.

The general setting is again a set of $s$ species, with $c_1, c_2, \dots, c_s$ representing their concentrations; a set of $n$ complexes, and a
 (numbered) set of reactions $E$ between these complexes. We also consider a finite directed graph $G= G(V,E,Y)$ whose vertices are labeled by the complexes and we keep our former notations.

We assume there is a non-negative continuous real-valued rate function $\mR_{ij}({\bc})$ for each reaction $(i,j)$ in the network, with the property that $\mR_{ij}({\bc})=0$ if and only if $c_k=0$ for any $k$ in the support of $y_i$. We record this information in the notation as $G=G(V,E,\mR,Y)$. In a mass-action kinetics system, we simply have
$\mR_{ij}({\bc}) = \kappa_{ij} {\bc}^{y_i}$, with $\kappa_{ij} \in \R_{>0}$, and in this case the notations
$G(V,E, \kappa, Y)$ and $G(V,E, \mR, Y)$ refer to the same system.

The instantaneous rate of change of $c_i$ is given by:
\begin{equation}\label{eq:geneq}
\frac{dc_i}{dt}=\underset{(j,i) \in E}{\sum}y_{ji}\mR_{ji}({\bc}) - \underset{(i,\ell) \in E}{\sum}y_{i\ell}\mR_{i\ell}({\bc}).
\end{equation}
The generalized system of differential equations that describe these dynamics could be written as
\begin{equation} \label{dynsysgeneral}
 \frac{d{\bc}}{dt}=\mR C_G^t Y,
\end{equation}
where $\mR$ is the $1 \times e$ matrix with entries $\mR_{ij}$, and $C_G^t$ is the
(transpose of) the corresponding signed incidence matrix we considered in subsection~\ref{ssec:la}.

\begin{remark}\label{rmk:notations}
It might be useful to compare our notation with the notation in \cite{karin,als07}.
First, we consider a row velocity vector $\frac{d{\bc}}{dt}$, while it is standard to consider the transposed column vector. So,
\[ \left(\frac{d{\bc}}{dt}\right)^t = \Gamma \mR({\bc})^t, \quad \text{where} \, \Gamma = Y^t.C_G \in \Z^{s \times e}.\]
Assume we have a mass-action kinetics system.
We denote by $K$ the $n \times e$ real matrix with entry in row indicated by complex $i$ and column indicated by the reaction edge $(i,j)$ equal to $\kappa_{ij}$, and equal to zero elsewhere. Then, $\mR({\bc}) = \Psi({\bc}) K$, the Laplace matrix equals $A_\kappa = K C_G^t$ and we have
\[\frac{d{\bc}}{dt} = \Psi({\bc}) A_\kappa Y = \Psi({\bc}) (K C_G^t) Y = \mR({\bc}) \Gamma^t.\]
In the notation of \cite{karin}, $Y^t$ is called $Y_s$, and the incidence matrices
 are denoted by $C_G = I_a, K^t = I_K$.
\end{remark}

In this general context we adapt the previous definitions.

\begin{definition}
 A \emph{complex balanced} kinetic system is a dynamical system ~\eqref{dynsysgeneral} associated with 
 the data $G(V,E, \mR, Y)$, for which the equations $\mR C_G^t = 0$ admit a strictly positive solution $\bc_0 \in \R_{>0}^s$ . Such a solution $\bc_0$ is a steady state of the system, i.e., the $s$ coordinates of $\mR C_G^t Y$ vanish. We call $\bc_0$ a complex balancing equilibrium.
\end{definition}

We will again assume that the digraph $G$ is reversible, and thus identify $G = (V, E)$ with the underlying undirected graph $\widetilde{G} = (V, \widetilde{E})$, where $\widetilde{E} = \{\{i, j\} : (i, j) \in E\}$.

\begin{definition}
A \emph{detailed balanced} kinetic system is a dynamical system ~\eqref{dynsysgeneral} asociated with the data $G(V,E, \mR, Y)$ for which the equations  $\mR_{ij}(\bc)-\mR_{ji}(\bc)= 0$, for all $\{i, j\} \in \widetilde{E}$, admit a strictly positive steady state $\bc_0 \in \R_{>0}^s$. We call $\bc_0$ a detailed balancing equilibrium.
\end{definition}

Again, every detailed balanced kinetic system is also complex balanced.
To define formal balancing, we need to start from a particular positive steady state:

\begin{definition} \label{def:formalgeneral}
Given a complex balanced system at the positive steady state $\bc_0 \in \R_{>0}^s$ corresponding to the
data $G(V,E, \mR, Y)$,
we say the kinetic system is \emph{formally balanced at $\bc_0$} (or that $\bc_0$ is a formally balancing
equilibrium) if the following condition holds for every cycle $\widetilde{C}$ of $\widetilde{G}$:
\begin{equation}\label{eq:fbqgeneral}
 \underset{(i,j) \; in \; C^+}{\prod}\mR_{ij}(\bc_0)=\underset{(j,i) \; in \; C^-}{\prod}\mR_{ji}(\bc_0).
\end{equation}
\end{definition}

We can now reformulate Theorem \ref{th:main}:

\begin{theorem} \label{th:maingeneral}
Consider a kinetic system ~\eqref{dynsysgeneral} asociated to the data $G(V,E, \mR, Y)$ with a complex balancing positive steady state $\bc_0 \in \R_{>0}^s$. We have that $\bc_0$ is a detailed balancing equilibrium  if and only if it the system is formally balanced at $\bc_0$.
\end{theorem}

\begin{proof}
 Given the complex balancing steady state $\bc_0 \in \R_{>0}^s$, we define constants $\kappa_{ij}=\mR_{ij}(\bc_0) \bc_0^{-y_i}$ for each $(i,j) \in E$ and we consider  the  mass-action kinetics dynamical system 
  $\frac{d{\bc}}{dt}=\Psi(\bc) A_\kappa Y$  associated with $G(V,E,\kappa,Y)$.
As $\mR_{ij}(\bc_0)=\kappa_{ij}\bc_0^{y_i}$, we have $\Psi(\bc_0)A_{\kappa}=0$, and so this new mass-action kinetics  system is complex balanced in the previous sense. 

Moreover, as the kinetic system is formally balanced at $\bc_0$, we have that
\[\underset{(i,j) \; in \; C^+}{\prod}\kappa_{ij} =  C \underset{(i,j) \; in \; C^+}{\prod}\mR_{ij}(\bc_0) =
C \underset{(j,i) \; in \; C^-}{\prod}\mR_{ji}(\bc_0)= \underset{(j,i) \; in \; C^-}{\prod}\kappa_{ji}, \]
where $C = c_0^{-\sum_{i \in E\left(\widetilde{C}\right)}y_i} \not =0$. Then, the mass-action kinetics system associated
with $G(V,E,\kappa,Y)$ is 
formally balanced. By Theorem \ref{th:main} it is detailed balanced. This means that every binomial $\kappa_{ij}\bc^{y_i}-\kappa_{ji}\bc^{y_j}$ vanishes at $\bc_0$, implying $\mR_{ij}(\bc)-\mR_{ji}(\bc_0)= 0$, and so the kinetic system associated with $G(V,E,\mR,Y)$ is detailed balanced at $\bc_0$. The other implication is clear.
\end{proof}

We end the paper by showing another necessary and sufficient condition for a complex balanced
kinetic system to be detailed balanced.

\begin{proposition}\label{prop:one}
 Given a kinetic system~\eqref{dynsysgeneral} asociated to the data $G(V,E, \mR, Y)$ with a complex balancing positive steady state $\bc_0 \in \R_{>0}^s$, the following statements are equivalent:
 \begin{itemize}
 \item[(i)] The equilibrium $\bc_0$ is detailed balancing.
 \item[(ii)] For every cycle $\widetilde{C}$ in $\widetilde{G}$ there exists an edge ${i_{\widetilde{C}},j_{\widetilde{C}}} \in E(\widetilde{C})$ such that $\mR_{ij}(\bc_0)-\mR_{ji}(\bc_0)= 0$. 
 \item[(iii)] Property (ii) holds for every \emph{basic} cycle associated to any spanning forest of $\widetilde{G}$.
 \end{itemize}
 \end{proposition}

\begin{proof}
The equivalence between (ii) and (iii) is clear, as well as the implication from
(i) to (ii). To see that (iii) implies (i), 
let $G'$ be the digraph obtained from $G$ by ``deleting'' all edges $(i_{\widetilde{C}},j_{\widetilde{C}}), {j_{\widetilde{C}},i_{\widetilde{C}}}$ in $E(C)$,
together with their corresponding labels, for all basic cycles $\widetilde{C}$.
Then, the associated undirected graph $\widetilde{G'}$ has no cycles and so 
any positive complex balancing  equilibrium $\bc_0$ for $G'$ is automatically also detailed balancing.
Call $A_\kappa(\bc_0)$ (resp. $A'_\kappa(\bc_0)$) the Laplace matrices of the mass-action kinetics system associated with $G$ (resp. $G'$) with reaction constants $\kappa_{ij}=\mR_{ij}(\bc_0)\bc_0^{-y_i} $ for each $(i,j) \in E$ (resp. $\kappa_{ij}=\mR_{ij}(\bc_0)\bc_0^{-y_i} $ for each $(i,j) \in E - \{(i_{\widetilde{C}},j_{\widetilde{C}}), (j_{\widetilde{C}},i_{\widetilde{C}}), \, \widetilde{C}
{\text{\, a \, basic \, cycle \, of \,}} \widetilde{G}\}$).
But if $\bc_0$ satisfies the conditions in (iii), it follows that
$$\Psi(\bc_0) A'_\kappa(\bc_0) \, = \, \Psi(\bc_0) A_\kappa(\bc_0) =0.$$
Therefore, $\bc_0$ is detailed balancing for $G'$, which together with the equalities
in (iii) implies that $\bc_0$ is detailed balancing for $G$, as wanted.
\end{proof}

\section{Conclusions}
We studied the conditions in parameter space which ensure the existence of particulary well behaved dynamics in general (mass-action) kinetics chemical reaction networks and we have clarified from an
algebraic perspective important classical notions. We plan to further apply this point of view to the study of biologically meaningful biochemical reaction networks, in particular those associated to enzymatic reactions as in \cite{gu07,gu08,als07,ctf06}, where our tools from elimination theory in the framework of algebraic varieties (and in particular, toric varieties), together from results in algebraic combinatorics as the Matrix-Tree Theorem allow to reformulate and generalize current approaches.

\bigskip

\noindent \textbf{{\large Acknowledgments:}}
We thank the Statistical and Applied Mathematical Sciences Institute (SAMSI), USA, where this work was initiated and completed. 
We are very grateful to Martin Feinberg for his thoughtful comments.

\medskip

\bibliographystyle{plain}

\end{document}